\documentclass[11pt,a4paper]{article}

\usepackage{amscd}       % Questo serve per fare i diagrammi
\usepackage{amsfonts}
\usepackage{amsmath}
\usepackage{amssymb}
\usepackage{amstext}
\usepackage{amsthm}       %proof environment :)
\usepackage{xspace}
\usepackage{graphicx}
\usepackage{url} 

\newcommand{\C}{\mathbb C}
\newcommand{\R}{\mathbb R}
\newcommand{\N}{\mathbb N}
\newcommand{\Q}{\mathbb Q}
\newcommand{\Z}{\mathbb Z} 

\newtheorem{theorem}{Theorem}[section]
\newtheorem{dfn}[theorem]{Definition}
\newtheorem{lemma}[theorem]{Lemma}
\newtheorem{cor}[theorem]{Corollary}
\newtheorem{rem}[theorem]{Remark}
\newtheorem{prop}[theorem]{Proposition}

\begin{document} 

\author{Paolo Ghiggini \footnote{The author is a member of EDGE, Research 
Training Network HPRN-CT-2000-00101, supported by The European Human Potential
 Programme. }\\
Mathematisches Institut, \\ Ludwig-Maximilians-Universit{\"a}t M{\"u}nchen \\
\url{ghiggini@mathematik.uni-muenchen.de}}
\title{Linear Legendrian curves in $T^3$} 
\date{}
\maketitle
%\tableofcontents
\begin{abstract}
Using convex surfaces and Kanda's classification theorem, we classify 
Legendrian isotopy classes of Legendrian linear curves in all tight contact 
structures on $T^3$. Some of the knot types considered in this article provide
new examples of non transversally simple knot types.  
\end{abstract}

\section{Introduction}
The study of Legendrian knots is an important topic in Contact Topology because
Legendrian knots are natural objects and capture the topology of the underlying
contact structure very well. In fact, for example, overtwistedness is defined 
in terms of Legendrian knots, and they have also proved useful in 
distinguishing tight contact structures. See, for example, 
\cite{ghiggini:1,honda:2,kanda:1}.

Recently, there have been some progress in the classification of Legendrian 
knots. Eliashberg and Fraser classified Legendrian unknots 
\cite{eliashberg-fraser} and,
using convex surfaces and bypasses, Etnyre and Honda classified Legendrian 
torus knots and figure-eight knots \cite{etnyre-honda:3}, and gave structure 
theorems for connected sums \cite{etnyre-honda:4} and cables  
\cite{etnyre-honda:5}.

In this article we classify Legendrian linear curves on $T^3$. This is the
first classification result for homologically non trivial Legendrian knots.
The contact structures $\xi_n= \ker \sin(2 \pi nz)dx+ \cos (2 \pi nz)dy$ on 
$\R^3 / \Z^3=T^3$
are all universally tight and pairwise non isomorphic by \cite{giroux:2}, and 
any tight contact structure on $T^3$ is isomorphic to one of them by 
\cite{kanda:1}. 
More specifically, any isotopy class of tight contact structures on $T^3$ has a 
representative obtained by the action of an element of $SL(3,\Z)$ on some
$\xi_n$. For this reason, classifying Legendrian linear curves in any tight
contact structure on $T^3$ is equivalent to classifying Legendrian linear 
curves in $(T^3, \xi_n)$ for any $n \in \N^+$.

The definition of the Thurston-Bennequin invariant and of the rotation number
can be extended to Legendrian linear knot in $T^3$, see \cite{kanda:1} and 
Section \ref{classici}, and provide necessary conditions for the existence of
a Legendrian isotopy. The flavour of the topology of Legendrian linear curves
in a given knot type is determined by the maximum value $\overline{\hbox{tb}}$
assumed by the 
Thurston-Bennequin invariant on the Legendrian curves in that knot type. 
If for a linear knot type 
$\overline{\hbox{tb}}<0$, then the 
the Thurston-Bennequin invariant and the rotation number are complete 
Legendrian isotopy invariants in that knot type. In the terminology of 
\cite{etnyre-honda:3} we say that the knot type is Legendrian simple. 

On the other hand, Legendrian linear curves in knot 
types with $\overline{\hbox{tb}}=0$ are more rigid.
The source of this dichotomy is the fact that a linear Legendrian curve with 
$tb<0$ can be embedded into a convex torus as a Legendrian ruling, while a 
linear Legendrian curve with $tb=0$ must be embedded as a Legendrian divide. 
A similar
phenomenon appeared in the classification of Legendrian torus knots 
\cite{etnyre-honda:3}.

By a theorem of Fuchs and Tabachnikov \cite{fuchs-tabachnikov}, any two 
smoothly isotopic Legendrian 
knots in $S^3$ are stably Legendrian isotopic, i.~e. they become 
Legendrian isotopic after some finite sequence of stabilisations. In general 
both positive and negative stabilisations are needed in the process, and
 their number can be arbitrarily high. In $T^3$ stable Legendrian Isotopy 
manifests itself with slightly different features. In fact any two Legendrian
linear curves with the same classical invariants become Legendrian isotopic 
after stabilising twice, once positively and once negatively, and some pairs of
Legendrian curves become Legendrian isotopic after stabilising both once in 
the same way. On the other hand, when $n>1$, some other pairs of linear 
Legendrian curves in $(T^3, \xi_n)$ with $\mbox{tb}=0$ remain non isotopic after
any number of stabilisations, if the stabilisations applied to each curve
have all the same sign. 
By \cite{etnyre-honda:3}, Theorem 2.10, this implies that the knot types with 
$\overline{\mbox{tb}}=0$ in $(T^3, \xi_n)$ are not transversally simple for $n>1$. 
Non  transversal simplicity is a recently 
discovered phenomenon, other examples have been previously found by Birman and 
Menasco \cite{birman-menasco}, and by Etnyre and Honda \cite{etnyre-honda:5}.

In this paper we follow the strategy proposed by Etnyre and Honda in 
\cite{etnyre-honda:3} in order to classify Legendrian knots. First, we 
classify Legendrian linear curves with 
maximal Thurston-Bennequin invariant. Then we prove that any Legendrian linear 
curve is a multiple stabilisation of a Legendrian linear curve with maximal 
Thurston-Bennequin
invariant, and finally we determine when Legendrian linear curves become
Legendrian isotopic under stabilisations.
 
\section{Statement of results}
Fix once and for all an identification $T^3 = (\R / \Z)^3$. The chosen 
identification induces global, $\Z$-periodic, coordinates $(x,y,z) \mapsto (y,z)$ 
on $T^3$. We introduce also a trivial bundle structure $T^3 \to T^2$ defined by the 
projection $(x,y,z)$. We will call an embedded torus (or annulus) 
{\em horizontal}, if it is isotopic to (a subset of) the image of a section, 
and {\em vertical} if it is incompressible and 
isotopic to a fibred torus (or annulus). This terminology will be used 
throughout the paper.
\begin{dfn}
For any $n \in \N_+$ we denote by $\xi_n$ the contact structure
$$\xi_n = \ker ( \sin (2 \pi nz)dx+ \cos (2 \pi nz)dy)$$
\end{dfn}
All these contact structures are universally tight and pairwise distinct
by \cite{giroux:2}.

\begin{dfn}
A simple closed curve $C \subset T^3$ is called {\em linear} if it is 
isotopic to a curve which lifts to a rational line in the universal cover
$\R^3 \to (\R / \Z)^3$. 
\end{dfn}
 The identification $T^3 = (\R / \Z)^3$ induces an identification 
$H_1(T^3, \Z)= \Z^3$ so that an oriented linear curve $C$ represents an homology 
class which corresponds to a primitive integer vector $(c_1,c_2,c_3)$. 
\begin{dfn}
 We will call the vector $(c_1,c_2,c_3)$ the {\em direction} of the linear curve $C$.
\end{dfn}
The direction
of a linear curve determines its isotopy type as an oriented curve completely.
\begin{dfn}
Given a primitive vector $(c_1,c_2,c_3) \in \Z^3$,
we define the knot type ${\cal K}(c_1,c_2,c_3)$ as the set of all linear curves 
with direction $(c_1,c_2,c_3)$.
\end{dfn}
\begin{dfn}
We define ${\cal L}_n(c_1,c_2,c_3)$ as the set of all linear curves 
$L \in {\cal K}(c_1,c_2,c_3)$ which are Legendrian with respect to the tight 
contact structure $\xi_n$.
\end{dfn}

Though linear curves in $T^3$ are not homologically trivial, the classical 
invariants for Legendrian knots in $S^3$ (Thurston-Bennequin invariant and 
rotation number) can be generalised 
to linear Legendrian curves in $(T^3, \xi_n)$, and give necessary conditions for
the existence of a Legendrian isotopy between smoothly isotopic Legendrian
linear curves. Their definitions and the proofs of their basic properties is
postponed to Section \ref{classici}. For these ``classical'' invariants
we will keep the standard terminology used for classical invariants in $S^3$.
Now we can state the main results of the article.

\begin{theorem} \label{main:tb<0} 
Let $L_1$ and $L_2$ be Legendrian curves in 
${\cal L}_n(c_1,c_2,c_3)$. If $c_3 \neq 0$, then $L_1$ is Legendrian isotopic to $L_2$ if 
and only if $\hbox{tb}(L_1)= \hbox{tb}(L_2)$ and $\hbox{r}(L_1)= \hbox{r}(L_2)$.
\end{theorem}

\begin{theorem}\label{main:tb=0}
Any linear Legendrian curve $L$ in ${\cal L}_n(c_1,c_2,0)$ has $\mbox{tb}(L) \leq 0$.
The set of the the Legendrian isotopy classes of Legendrian curves in 
${\cal L}_n(c_1,c_2,0)$ with $\hbox{tb}=0$ is in bijection with the set $\pi_0(\Gamma_T)$ 
of the dividing curves on a horizontal convex $\# \Gamma$-minimising torus $T \subset T^3$.
Any Legendrian curve $L \in {\cal L}_n(c_1,c_2,0)$ with $\hbox{tb}(L)<0$ is 
Legendrian 
isotopic to a (multiple) stabilisation of a Legendrian curve $L'$ with 
$\hbox{tb}(L')=0$.
Two Legendrian curves $L_1$ and $L_2$ in ${\cal L}_n(c_1,c_2,0)$ are Legendrian 
isotopic if and only if $\hbox{tb}(L_1)= \hbox{tb}(L_2)$, $\hbox{r}(L_1)= 
\hbox{r}(L_2)$ and at least one of the followings holds.
\begin{enumerate}
\item $L_1'$ is Legendrian isotopic to $L_2'$,
\item $|\hbox{r}(L_i)| < |\hbox{tb}(L_i)|$,
\item $|\hbox{r}(L_i)| = |\hbox{tb}(L_i)|$, and $L_1'$ and $L_2'$ correspond to
dividing curves of $T$ bounding a region $C$, and the sign of 
$C$ is opposite to the sign of $\hbox{r}(L_i)$.
\end{enumerate}
\end{theorem}

\section{``Classical'' invariants} \label{classici}
The definition of the classical invariants for the linear Legendrian curves in 
$T^3$ is somehow common folklore among the Contact Topology people, however, a 
careful exposition of the topic is hard to find in the literature. For this 
reason, we have decided to include the proofs of the principal results.

\begin{dfn} \label{tb} (\cite{kanda:1}, Definitions 7.2 and 7.3).
Let $L \subset (T^3, \xi_n)$ be a linear Legendrian curve, and let $T \subset T^3$ be an 
incompressible torus containing $L$. Then the {\em Thurston-Bennequin invariant}
$\hbox{tb}(L)$ is defined as the 
twisting of the framing of $L$ induced by $\xi_n$ with respect 
to the framing induced by the torus $T$.
\end{dfn}
{\em A priori}, the Thurston-Bennequin Invariant of a Legendrian linear curve
 could depend on the isotopy class of the torus
$T$, but it turns out that this is not the case. 

\begin{lemma}\label{framing} 
The Thurston-Bennequin Invariant of a linear Legendrian curve $\hbox{tb}(L)$ 
does not depend on the torus used to compute it.
\end{lemma}

\begin{proof}
If $L$ is a linear curve, then there is a projection $\pi: T^3 \to T^2$ such that 
$L$ is a fibre of $\pi$. We will show that, for any incompressible torus $T$ 
containing $L$,
the framing on $L$ induced by $T$ is homotopic to the product framing induced 
by $\pi$.
 
 Let $\R^2 \to T^2$ be the universal cover, and let $\tilde{\pi}: \R^2 \times S^1 \to \R^2$ 
be the pull-back of $\pi$. The torus $T$ lifts 
to an annulus $A$ with boundary $\partial A= \tilde{L}_0-\tilde{L}_1$, where $\tilde{L}_0$
and $\tilde{L}_1$ are connected components of the pre-image of $L$.
The difference between the framing induced by $T$
 and the projection framing is the number of intersection points, counted with
 sign, between $A$ and a fibre 
$\tilde{L}'$ of $\tilde{\pi}$ near to 
$\tilde{L}_0$. But every fibre is homologous to a fibre disjoint from $A$,
therefore the difference between the two framings is zero. 
\end{proof}

\begin{dfn}
Let $(c_1,c_2,c_3)$ be a primitive integer vector.  We define the 
{\em maximum Thurston-Bennequin invariant} 
$\overline{\hbox{tb}}({\cal L}_n(c_1,c_2,c_3))$ of the 
knot type ${\cal K}(c_1, c_2, c_3)$ for the tight contact structure $\xi_n$  as
$$\overline{\hbox{tb}}({\cal L}_n(c_1,c_2,c_3))= \sup \{ tb(L): L \in {\cal L}_n(c_1,c_2,c_3) \}$$
\end{dfn} 
In general, the maximum Thurston-Bennequin invariant in a given knot type is 
not easy to 
compute, however, in the case of linear curve, it can be expressed  in 
terms of the direction of the curve and of the contact structure.

\begin{theorem}\label{kanda}(\cite{kanda:1}, Theorem 7.6).
The maximum Thurston-Bennequin invariant of the knot type ${\cal K}(c_1,c_2,c_3)$
for the contact structure $\xi_n$ is 
$$\overline{\hbox{tb}}({\cal L}_n(c_1,c_2,c_3)) = -|nc_3|$$
\end{theorem}  

\begin{dfn}
Let $X$ be a nowhere vanishing section of $\xi_n$, and let $L \subset (T^3, \xi_n)$ be an 
oriented Legendrian knot. We define the
{\em rotation number} of $L$ with respect to $X$, denoted by $\hbox{r}_X(L)$, 
as the twisting of the tangent vector $\dot{L}$ with respect to $X$ along $L$.
\end{dfn}  

The definition of the rotation number makes sense because $c_1(\xi_n)=0$ implies
 the existence of 
nowhere vanishing sections of $\xi_n$. However, the rotation number does 
depend on the choice of the vector field $X$ as follows.
Let $X$ and $Y$ be nowhere vanishing sections of $\xi_n$. The choice of a complex
structure on $\xi_n$ determines a function $u:T^3 \to \C^*$ such that $Y_p=u(p)X_p$ for 
any $p \in T^3$. Let $c(X,Y) \in H^1(T^3, \Z)=[T^3, S^1]$ be the cohomology class 
corresponding
to the homotopy class of $u$, then $\hbox{r}_Y(L)= \hbox{r}_X(L)+ \langle c(X,Y),[L] \rangle$.

\begin{dfn}
Let $L$ and $L'$ be isotopic Legendrian knots in $T^3$, and let $\Sigma \subset T^3$ be
an embedded surface with boundary $\partial \Sigma = L - L'$. We define the 
{\em relative rotation number} of $L$ with respect to $L'$ $\hbox{r}(L,L')$ as the 
obstruction to extending $\dot{L}$ and $\dot{L}'$ to a 
trivialisation of $\xi_n|_{\Sigma}$.
\end{dfn}

The relative rotation number does not depend on the surface bounding $L$ and 
$L'$ because $c_1(\xi_n)=0$. Moreover,  for any nowhere 
vanishing section $X$ of $\xi_n$, the equality
$$\hbox{r}(L, L')= \hbox{r}_X(L)- \hbox{r}_X(L')$$
holds. 

\begin{lemma} \label{importante}
If $L$ and $L'$ are smoothly isotopic linear Legendrian curves in $(T^3, \xi_n)$
realising the maximum of the Thurston-Bennequin invariant in their knot type,
 then $\mbox{r}(L,L')=0$.
\end{lemma}

\begin{proof}
Let $\Sigma \subset T^3$ be a convex embedded surface bounding $L$ and $L'$. 
First, we observe that the equality $\mbox{r}(L,L')=0$ holds if $\Sigma$ is an 
annulus. In fact $\Gamma_{\Sigma}$ contains no boundary parallel dividing arc because
$\mbox{tb}(L)= \mbox{tb}(L')$ is maximal in the knot type of $L$ 
and $L'$, thus $\mbox{r}_{\Sigma}(L,L')= \chi(\Sigma_+)- \chi(\Sigma_-)=0$. If $\Sigma$ is not an annulus, 
take a finite
cover $\pi: T^3 \to T^3$ which is trivial when restricted to $L$ and $L'$, and choose
 a lift $\widehat{\Sigma}$ of $\Sigma$. Let $\widehat{\xi}$ the pulled-back tight contact 
structure, and $\partial \widehat{\Sigma}= \widehat{L} \cup \widehat{L}'$. If the 
cover $\pi: T^3 \to T^3$ has sufficiently many sheets, we can assume that there
exists a $\mbox{tb}$-maximising Legendrian linear curve $\widehat{L}_0$ such 
that $\widehat{L}$ and $\widehat{L}_0$ bound an annulus $A$ and $\widehat{L}'$ 
and $\widehat{L}_0$ bound an annulus $A'$.   Then 
$$r_{\widehat{\Sigma}}(\widehat{L}, \widehat{L}')+ r_{A'}(\widehat{L}', \widehat{L}_0)+
r_A(\widehat{L}_0, \widehat{L})= \langle c_1(\widehat{\xi}), [\widehat{\Sigma} \cup A \cup A'] \rangle =0$$
and this implies $r_{\Sigma}(L, L')=0$ because $r_{A'}(\widehat{L}', \widehat{L}_0)= 
r_A(\widehat{L}_0, \widehat{L})= 0$.
\end{proof}

\begin{rem}
Lemma \ref{importante} allows us to choose the section $X$ so that any
$\mbox{tb}$-maximising Legendrian linear curve $L$ in a given knot type has
rotation number $\mbox{r}_X(L)=0$. In the following, whenever this choice is 
done, we will write simply $\mbox{r}(L)$, omitting the section. 
\end{rem}

\begin{figure}\centering
\includegraphics[width=5cm]{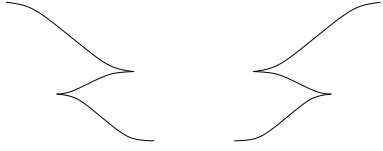}
\caption{Zig-zags in the front projection.}
\label{zigzag.fig}
\end{figure}
Given a Legendrian knot $L$, there is a standard procedure called 
{\em stabilisation} which produces another Legendrian knot in the same knot type.
In $\R^3$ with the standard contact form $\alpha_0= dz-ydx$ any Legendrian knot can be
recovered from its projection on the $xz$--plane (the front projection) by the 
equation $y=\frac{dz}{dx}$. Stabilisation is performed in the standard $\R^3$ 
by adding a zig-zag to the front projection of the Legendrian knot, and the 
definition can be extended to Legendrian knots in any contact manifold by 
Darboux Theorem.
See for example \cite{etnyre:0}, Section 2.7, or \cite{etnyre-honda:3}, 
Section  2.1. A stabilisation comes with a sign, which corresponds to the 
choice of one of the two possible ways to add a zig-zag to the front 
projection. See Figure \ref{zigzag.fig}. We denote by $S_+$ a positive
stabilisation, and by $S_-$ a negative one. Positive and negative stabilisation
are well defined up to Legendrian isotopy, and commute.  
See \cite{etnyre-honda:3}, Lemma 2.5.
Stabilisations act on the classical invariants in the following way:
\begin{align*}
\mbox{tb}(S_{\pm}(L)) &= \mbox{tb}(L)-1 \\
\mbox{r}(S_{\pm}(L)) &= \mbox{r}(L) \pm 1
\end{align*}

In the following, we will call a composition of stabilisations a 
{\em multiple stabilisation} and sometimes, when we want stress the point that
we are considering only one stabilisation, we will call a stabilisation
a {\em simple stabilisation}.

The converse of a stabilisation is called {\em destabilisation}. Contrarily to 
stabilisations, destabilisations are not always possible. In fact a 
destabilisation is equivalent to finding a bypass: see \cite{etnyre:0}, Lemma
2.20. 

%%%%%%%%%%%%%%%%%%%%%%%%%%%%%%%%%%%%%%%%%%%%%%%%%%%%%%%%%%%%%%%%%%%%%%%%%%%%%%
\section{Legendrian curves with $\overline{\hbox{tb}}<0$}
In this section we prove Theorem \ref{main:tb<0}.  
The proof will proceed in two steps. First,
we will show that all Legendrian curves realising the maximum twisting number 
in a given linear knot type with $\overline{\hbox{tb}}<0$ are Legendrian
 isotopic. Then, we will show that any Legendrian curve in that  linear knot 
type is Legendrian isotopic to a (multiple) stabilisation of a Legendrian 
curve with maximum twisting number. For the rest of the section, 
${\cal K}(c_1,c_2,c_3)$ will denote a linear knot type with $c_3 \neq 0$, so that  
$\overline{\hbox{tb}}({\cal L}_n(c_1,c_2,c_3))= - n|c_3|<0$.

Controlling isotopies  between Legendrian knots is much harder than 
controlling isotopies between incompressible tori containing them. The 
following two lemmas describe the dividing set of convex incompressible tori 
in $(T^3, \xi_n)$.

\begin{lemma}\label{divisione}
Let $\left ( \begin{array}{c}
             b_1 \\ b_2 \\ b_3
 \end{array} \right )$ and $\left ( \begin{array}{c}
             c_1 \\ c_2 \\ c_3
 \end{array} \right )$
be two vectors of an integer basis of $\Z^3$, and
 let $T \subset (T^3, \xi_n)$ be the image of the embedding $\iota:  T^2 \to T^3$ 
induced by 
$$\tilde{\iota} : \R^2 \ni \left ( \begin{array}{l} u \\ v \end{array}
\right ) \mapsto \left ( \begin{array}{l} 
ub_1+vc_1 \\ ub_2+vc_2 \\ ub_3+vc_3 \end{array} \right ) \in \R^3$$ 

Then $T$ is a convex torus with slope $s(T)=- \frac{b_3}{c_3}$ and 
$\# \Gamma_T=2n(b_3,c_3)$, where
$(b_3,c_3)$ is the greatest common divisor between $b_3$ and $c_3$. 
\end{lemma}
\begin{proof}
  Since $c_3 \neq 0$, there exist constants $a$, $b \in \R$ such that 
$X= a \frac{\partial}{\partial x} + b \frac{\partial}{\partial y}$ is transverse to $T$. For any $n$, the 
contact structures 
$\xi_n$ are invariants in the directions $x$ and $y$, therefore $T$ is convex
 because  $X$ is a contact vector field.  
 The characteristic set $\Sigma = \{ p \in T^3 \; | \; X(p) \in \zeta_n(p) \}$ consists of  $2n$ 
parallel copies of a horizontal torus of the form $\{z \in \Z \}$, therefore
 the dividing set $\Gamma_T= T \cap \Sigma$ of $T$ consists of $2n$ copies of the
the set $\{ub_3+vc_3 \in \Z \}$. By an easy computation, $s(T)=- \frac{b_3}{c_3}$ and 
$\# \Gamma_T=2n(b_3,c_3)$.
\end{proof}
\begin{lemma}\label{preparazione}
Let $T$ be a vertical torus in standard form as in 
Lemma \ref{divisione} with Legendrian rulings with direction $(c_1, c_2, c_3)$, 
and let $S$, $S'$ be smoothly isotopic convex tori 
intersecting $T$ along a ruling of $T$. Then $S$ is contact isotopic to $S'$.
\end{lemma}
\begin{proof}
If $T \cap S$ and $T \cap S'$ are distinct ruling curves of $T$, there is a vector
field tangent to $T$ whose flow maps $T \cap S$ onto $T \cap S'$, and which preserves
the characteristic foliation of $T$. This vector field can be extended to a 
contact vector field on $T^3$ whose flow maps $S'$ onto a torus $S''$ such that
$T \cap S''= T \cap S$, therefore we can assume without loss of generality that 
$T \cap S'= T \cap S$. 

By the Isotopy Discretisation Lemma \cite{gluing}, Lemma 3.10, there is a 
sequence of standard tori $S=S_0, \ldots ,S_m=S'$ such that $S_{i-1} \cap T= S_i \cap T$ for 
$i=1, \ldots ,m$, $S_{i-1}$ is disjoint from $S_i$ outside $T$ and $S_i$ is obtained from
$S_{i-1}$ by a single bypass attachment. We will prove that $\Gamma_{S_{i-1}}$ is 
isotopic to $\Gamma_{S_i}$, then the bypass attachment from $S_{i-1}$ to $S_i$ is trivial
and, by \cite{gluing} Lemma 2.10, $S_{i-1}$ is contact isotopic to $S_i$, for 
$i= 1, \ldots ,m$. In the following we will consider $(T^3 \setminus T, \xi_n|_{T^3 \setminus T})$ as a 
nonrotative thickened torus $(T^2 \times I, \xi)$ and identify the tori $S_i$ with 
convex vertical annuli with Legendrian boundary in $(T^2 \times I, \xi)$.

Passing to a finite covering of $(T^2 \times I, \xi)$ in the horizontal direction, and 
lifting $S_i$, for $i=1, \ldots ,n$ to the covering if necessary, we can always 
assume that $(T^2 \times I, \xi)$ has integer boundary slopes.
If $(T^2 \times I, \xi)$ has integer boundary slopes, the isotopy between $\Gamma_{S_{i-1}}$
and $\Gamma_{S_i}$ is a consequence of the classification of non rotative tight
contact structures on thickened tori \cite{honda:1} Lemma 5.7. 
\end{proof}

\begin{lemma}\label{minimale}
Let $T$ and $S$ be non isotopic convex incompressible tori in standard form
intersecting along a common Legendrian ruling with direction $(c_1, c_2, c_3)$ with
$c_3 \neq 0$. Then, for any pair of convex tori $T'$ and $S'$ intersecting 
along a common vertical Legendrian ruling of $T'$ and isotopic 
to $T$ and $S$ respectively, 
\begin{enumerate}
\item $s(T')=s(T)$ and $\# \Gamma_{T'} \geq \# \Gamma_{T}$,
\item $S'$ is contact isotopic to $S$.
\end{enumerate}
In particular, $S'$ is $\# \Gamma$-minimising.
\end{lemma}
\begin{proof}
By Isotopy Discretisation \cite{gluing}, Lemma 3.10, there
is a sequence of convex tori $T=T_0, \ldots ,T_n=T'$ in standard form with Legendrian
rulings with direction $(c_1, c_2, c_3)$ such that $T_i$ is 
obtained from $T_{i-1}$ by attaching a bypass. In particular, $T^3 \setminus (T_{i-1} \cup T_i)=
N_i' \cup N_i''$, where $N_i$ and $N_i''$ are both diffeomorphic to $T^2 \times I$. We can 
assume inductively that $T_{i-1}$ satisfies the following assumptions:
\begin{enumerate}
\item $s(T_{i-1})=s(T)$ and $\# \Gamma_{T_{i-1}} \geq \# \Gamma_{T}$,
\item there is a convex vertical annulus $A_{i-1} \subset M \setminus T_{i-1}$ with Legendrian 
boundary on $\partial (T^3 \setminus T_{i-1})$ such that $A_{i-1}$ closes to a convex 
torus $S_{i-1} \subset M$ contact isotopic to $S$.
\end{enumerate}

A priori there are two kinds of transitions from $T_{i-1}$ to $T_i$:
\begin{enumerate}
\item $\# \Gamma_{T_i} = \# \Gamma_{T_{i-1}}= 2$ and $s(T_i) \neq s(T_{i-1})$
\item  $s(T_i) = s(T_{i-1})$ and $\# \Gamma_{T_i} = \# \Gamma_{T_{i-1}} \pm 2$
\end{enumerate} 

{\bf Case 1.} We will show that there are no transitions which change the 
slope. Suppose by contradiction that $\# \Gamma_{T_{i-1}}= \# \Gamma_{T_i}=2$ and 
$s(T_i) \neq s(T_{i-1})$, then $T^3 \setminus T_{i-1}$ is contactomorphic to a rotative thickened
torus. By \cite{honda:1}, Proposition 4.16, there is a convex torus parallel to
$T_{i-1}$ with slope $s$ for any $s \in \Q \cup \{ \infty \}$. In particular, there is a 
standard torus whose Legendrian divides have direction $(c_1, c_2, c_3)$. This 
contradicts Theorem \ref{kanda}, because $\overline{\mbox{tb}}({\cal L}_n(c_1, c_2,
c_3))=-n|c_3| < 0$.

{\bf Case 2.}
We will show that no
transition can decrease $\# \Gamma_{T_{i-1}}$ below $\# \Gamma_{T}$ first. Suppose that 
$\# \Gamma_{T_i} < \# \Gamma_{T_i-1} = \# \Gamma_{T_0}$, and let $L_0$, $L_i$ be  Legendrian rulings 
 of $T_0$ and $T_i$ respectively, both with direction $(c_1, c_2, c_3)$. Since 
$s(T_i)=s(T_0)$ and $\# \Gamma_{T_i}< \# \Gamma_{T_0}$, we have 
$$\mbox{tb}(L_i) = - \frac 12 |L_i \cap \Gamma_{T_i}| > - \frac 12 | L_0 \cap \Gamma_{T_0} | = 
\mbox{tb}(L_0) = -n |c_3|$$ contradicting Theorem \ref{kanda}. 

Take convex annuli $A_i' \subset N_i'$ and $A_i'' \subset N_i''$ with 
boundary on the same Legendrian rulings of $T_{i-1}$ and $T_i$, and define 
$S_i= A_i' \cup A_i'' \subset T^3$. 
By Proposition \ref{preparazione}, $S_i$ is contact isotopic to $S_{i-1}$ because 
both intersect $T_{i-1}$ along a Legendrian ruling. 
\end{proof} 

\begin{lemma}\label{torisotopi}
Let $T$ and $T'$ be smoothly isotopic $\# \Gamma$-minimising standard tori
with Legendrian rulings with direction $(c_1, c_2, c_3)$ with $c_3 \neq 0$. Then $T$ 
is contact isotopic to $T'$.  
\end{lemma}
\begin{proof}
Let $S$ and $S'$ be standard tori intersecting $T$ and $T'$ along a 
Legendrian ruling of $T$ and $T'$ respectively, as required by Lemma 
\ref{minimale}, then $S$ and $S'$ are contact isotopic. 
If there is a boundary parallel dividing arc in $S \setminus T$, we can perturb the 
characteristic foliation of $S$ in order to obtain a bypass attached 
vertically to $T$. 
This bypass would either change the slope of $T$, contradicting \ref{minimale},
or decrease the division number of $T$, contradicting the minimality 
assumption. Then we can make $S$ standard with vertical ruling, so that $S \cap T$ is a common
vertical ruling for both tori. Since the same argument applies to $T'$ and $S'$
too there is a complete symmetry between $T$ and $S$, then from Lemma 
\ref{minimale} we conclude that $T'$ is isotopic to $T$.
\end{proof}

\begin{prop}\label{massimali:tb<0}
  Let $L_1$ and  $L_2$  be linear Legendrian curves in ${\cal L}_n(c_1,c_2,c_3)$ with 
$\hbox{tb}(L_1)= \hbox{tb}(L_2)= \overline{\hbox{tb}}({\cal L}_n(c_1,c_2,c_3))<0$, 
then $L_1$ and $L_2$ are Legendrian isotopic.
\end{prop}
\begin{proof}
Let $T_1$ and $T_2$ be isotopic incompressible convex tori containing $L_1$ and 
$L_2$ respectively, then 
$$L_1 \cap \Gamma_{T_1} = \hbox{tb}(L_1)= \hbox{tb}(L_2)= L_2 \cap \Gamma_{T_2}$$
 For $i=1,2$, $L_i$ must intersect $\Gamma_{T_i}$ minimally, 
and both $T_1$ and $T_2$ are $\# \Gamma$-minimising in their isotopy class because $L_1$ 
and $L_2$ are $\hbox{tb}$-maximising in their knot type. Put $T_1$ and $T_2$ in 
standard form so that $L_1$ and $L_2$ are Legendrian rulings. By Lemma
\ref{torisotopi}, there is a contact isotopy 
$\varphi_t: (T^3, \xi_n) \to (T^3,\xi_n)$ such that $\varphi_0= \hbox{id}$ and $\varphi_1(T_1)=T_2$. The curve
$\varphi_1(L_1)$ is a Legendrian ruling of $T_2$, therefore it is Legendrian isotopic
to $L_1$ through Legendrian rulings of $T_2$.
\end{proof}

\begin{prop}\label{stabilizzazione:t<0}
Let $L \in  {\cal L}_n(c_1,c_2,c_3)$ be a linear Legendrian curve with
$\hbox{tb}(L) < \overline{\hbox{tb}}({\cal L}_n(c_1,c_2,c_3))$. Then there is 
a Legendrian linear curve $L' \in  {\cal L}_n(c_1,c_2,c_3)$ with $\hbox{tb}(L')= 
\overline{\hbox{tb}}({\cal L}_n(c_1,c_2,c_3))$ such that $L$ is a stabilisation of
$L'$.
\end{prop}
\begin{proof}
Let $T \subset (T^3, \xi_n)$ be an incompressible convex torus  containing 
$L$, and take a translate torus $T'$ contained in an invariant neighbourhood of
$T$. Put $T'$ in standard form with rulings with direction $(c_1, c_2, c_3)$, and 
consider a convex annulus $A$ between $L$ and a  
Legendrian ruling $R$ of $T'$. If $A$ is completely contained
in an invariant neighbourhood of $T$, then $\Gamma_A$ contains no
dividing arcs with both endpoints on $T'$, therefore $L$ is Legendrian 
isotopic to a stabilisation 
of $R$. Take now a convex torus $S$ intersecting $T'$ along $R$, and a 
translate $S'$ of $S$ contained in an invariant neighbourhood of $S$. As 
before, make $S'$ standard with Legendrian rulings with direction 
$(c_1, c_2, c_3)$, then $R$ is 
Legendrian isotopic to a stabilisation of a Legendrian ruling $L'$ of $S'$. 
Since $S'$ is $\# \Gamma$-minimising by Lemma \ref{minimale}, $L'$ realises the 
maximal twisting number in ${\cal L}_n(c_1, c_2, c_3)$.
\end{proof} 
\noindent
{\em Proof of Theorem \ref{main:tb<0}.}
Let $L_1$ and $L_2$ be Legendrian curves in ${\cal L}_n(c_1,c_2,c_3)$ such that
$\hbox{tb}(L_1)= \hbox{tb}(L_2)$ and $\hbox{r}(L_1)= \hbox{r}(L_1)$. Then, by 
Proposition \ref{stabilizzazione:t<0}, there are Legendrian curves
$L_1'$ and $L_2'$ realizing the maximum of the Thurston-Bennequin invariant, 
such that $L_i$ is a (multiple) stabilisation of $L_i'$, for $i=1,2$. By 
Proposition \ref{massimali:tb<0}, $L_1'$ is Legendrian isotopic to $L_2'$ by 
an isotopy induced by a contact isotopy $\varphi_t$ of $(T^3, \xi_n)$. The 
curves $L_2$ and $\varphi_t(L_1)$ are both (multiple) stabilisations of $L_2'$, then they 
are Legendrian isotopic because they have the same classical invariants. 
\qed \\

%%%%%%%%%%%%%%%%%%%%%%%%%%%%%%%%%%%%%%%%%%%%%%%%%%%%%%%%%%%%%%%%%%%%%%%%%%%%
\section{Legendrian curves with $\overline{\hbox{tb}}=0$}
In this section we classify the Legendrian linear curves in 
${\cal L}_n(c_1,c_2,0)$. First, we classify linear the Legendrian curves with 
$\mbox{tb}=0$, then we prove that any linear Legendrian curve
in ${\cal L}_n(c_1,c_2,0)$ is a stabilisation of a linear Legendrian curve with
$\mbox{tb}=0$. The final part of the section is devoted to prove the necessary
and sufficient conditions for two stabilisations of non Legendrian isotopic
linear Legendrian curves to become Legendrian isotopic.

By \cite{kanda:1}, Proposition 7.5, there is a matrix $A \in
SL(3, \Z)$ such that $A({\cal K}(c_1,c_2,0))= {\cal K}(1,0,0)$ and $A_*(\xi_n)$ is
isotopic to $\xi_n$, therefore the classification of the linear Legendrian  
curves in ${\cal L}_n(1,0,0)$ is equivalent to the classification of the linear 
Legendrian curves in ${\cal L}_n(c_1,c_2,0)$ for any primitive vector 
$(c_1,c_2) \in \Z^2$.
 Linear curves in ${\cal K}(1,0,0)$ will be called {\em vertical curves}. 

The contact structures $\xi_n$ are $S^1$-invariant with respect to the $S^1$-bundle 
structure on $T^3$ defined by the projection $(x,y,z) \mapsto (y,z)$. let $\Sigma  \subset T^3$ be the 
horizontal torus $\Sigma = \{ x=0 \}$, then $\Sigma$ is convex and $\Gamma_{\Sigma}$ consists of 
$2n$ closed curves. We will call $\Gamma_{T^2}$ the image of $\Gamma_{\Sigma}$ under the 
projection $T^3 \to T^2$.
For any $p \in \Gamma_{T^2}$  the curve  $L_0= S^1 \times \{ p \}$ is a Legendrian fibre
with twisting number $\hbox{tb}(L_0)=0$.

\begin{prop} \label{massimali-distinzione:t=0}
Let $p_0$, $p_1 \in \Gamma_{T^2}$, then $S^1 \times \{ p_0 \}$ is Legendrian isotopic to $S^1 \times 
\{ p_1 \}$ if and only if $p_0$ and $p_1$ belong to the same connected component of 
$\Gamma_{T^2}$.
\end{prop}
\begin{proof}
If $p_0$ and $p_1$ belong to the same connected component of $\Gamma_{T^2}$, there is
a path $\{ p_t \}_{t \in [0,1]} \subset \Gamma_{T^2}$ joining $p_0$ and $p_1$, then the family of 
Legendrian fibres 
$S^1 \times \{ p_t \}$ provides a Legendrian isotopy between $S^1 \times \{ p_0 \}$ and $S^1 \times 
\{ p_1 \}$. 

On the other hand, assume by contradiction that $S^1 \times \{ p_0 \}$ and $S^1 \times \{ p_1 \}$ 
are
Legendrian isotopic, but $p_0$ and $p_1$ belong to different connected components
of $\Gamma_{T^2}$. After possibly lifting the Legendrian isotopy between $S^1 \times \{ p_0 \}$
and $S^1 \times \{ p_1 \}$ to a suitable finite covering $(T^3, \widehat{\xi})$ of $(T^3 \xi_n)$,
we can assume that 
this Legendrian isotopy is disjoint from $S^1 \times \{ p_2 \}$, for some $p_2$ in the 
same dividing curve as $p_0$.
We can extend the Legendrian isotopy 
between $S^1 \times \{ p_0 \}$ and $S^1 \times \{ p_1 \}$ to a contact isotopy $\varphi_t$ fixed in
neighbourhood of $S^1 \times \{ p_2 \}$.
 For $i=0,1$, let $M_i= T^3 \setminus (\nu(S^1 \times \{ p_i \}) \cup \nu(S^1 \times \{ p_2 \}))$, where 
$\nu(S^1 \times \{ p_i \})$ is a 
standard neighbourhood of $S^1 \times \{ p_i \}$, then $\varphi_1$ is a contactomorphism between
$(M_0, \widehat{\xi}|_{M_0})$ and $(M_1, \widehat{\xi}|_{M_1})$.  This contradicts the 
classification
of $S^1$-invariants tight contact structures, \cite{honda:2}, Proposition 4.4, 
because the dividing sets $\Gamma_{\Sigma_0}$ and $\Gamma_{\Sigma_1}$ of convex $\# \Gamma$-minimising 
sections $\Sigma_0$ and $\Sigma_1$ of $(M_0, \widehat{\xi}|_{M_0})$ and $(M_1, \widehat{\xi}|_{M_1})$
respectively, obtained by removing small discs from $T^2$ around $p_i$ and $p_2$, 
for $i=0,1$, cannot be diffeomorphic.
\end{proof}

\begin{prop} \label{massimali-modelli:t=0}
If $L \subset (T^3, \xi_n)$ is a vertical Legendrian curve with $\hbox{tb}(L)=0$, then
there is a point $p \in \Gamma_{T^2}$ such that $L$ is Legendrian isotopic to $S^1 \times 
\{ p \}$.
\end{prop}
\begin{proof}
Let $T \subset T^3$ be the fibred torus $T = \{ y=0 \}$, and let $T_0$ be 
a convex vertical torus in standard form isotopic to $T$ containing $L$ as a 
Legendrian divide. By isotopy discretisation, there is a family of standard
tori $T_0, \ldots ,T_n=T$ such that $T_{i-1} \cap T_i = \emptyset$ for $i=1, \ldots ,n$.
For any $i=0, \ldots ,n$, the torus $T_i$ can be put in 
standard form with infinite slope and horizontal rulings because the projection
of $T$ intersect $\Gamma_{T^2}$, and $T^3 \setminus T_i$ is non rotative for any $i$.
In order to prove the proposition, we prove inductively that any Legendrian 
divide of $T_{i-1}$ is isotopic to a Legendrian divide of $T_i$.

Let $N_i', N_i'' \subset T^3$ be the connected components of $T^3 \setminus (T_{i-1} \cup T_i)$, and let 
$\Sigma_i \subset T^3$ be convex $\# \Gamma$-minimising horizontal tori such that $T_{i-1} \cap \Sigma_i$ and 
$T_i \cap \Sigma_i$ are Legendrian rulings of $T_{i-1}$ and $T_i$. By \cite{honda:2}, 
Proposition 4.4, $\Gamma_{\Sigma_i}$ is isotopic to $\Gamma_{\Sigma_0}$ for any $i=0, \ldots ,n'$, and by
\cite{honda:1}, Lemma 5.7, for any 
$i=1, \ldots ,n$, $\xi_n|_{N_i'}$ is isomorphic to an $S^1$-invariant contact structure 
on $S^1 \times (\Sigma_i \cap N_i')$ and $\xi_n|_{N_i''}$ is isomorphic to an $S^1$-invariant contact 
structure on $S^1 \times (\Sigma_i \cap N_i'')$. Any Legendrian divide of $T_{i-1}$ corresponds to
a point $p_{i-1} \in T_{i-1} \cap \Gamma_{\Sigma_i}$, either dividing arc of $\Gamma_i$ starting from 
$p_{i-1}$ eventually hits $T_i$ because $|T_i \cap \Sigma_i|=|T_0 \cap \Sigma_0| \neq 0$, where 
$| \cdot |$ denotes the minimal geometric intersection. Let $t \mapsto p_{i-1+t}$ be a 
parametrisation of a segment of dividing arc starting from $p_{i-1}$ such that 
$p_i \in T_i \cap \Gamma_{\Sigma_i}$. The path $t \mapsto p_{i-1+t}$ provides the required isotopy between 
a Legendrian divide of $T_{i-1}$ and a Legendrian divide of $T_i$.
\end{proof}

\begin{cor}
The set of the Legendrian isotopy classes of Legendrian curves $L$ in 
${\cal L}_n(1,0,0)$ with $\hbox{tb}(L)=0$ is  in bijection with the set 
$\pi_0(\Gamma_{T^2})$ of the connected components of $\Gamma_{T^2}$.
\end{cor}

\begin{prop}\label{stabilizzazione:t=0}
If $L \subset (T^3, \xi_n)$ is a vertical Legendrian curve with $\hbox{tb}(L)< 0$, then
$L$ is isotopic to a stabilisation of a vertical Legendrian curve $L_0$ with
$\hbox{tb}(L_0)= 0$.
\end{prop}

\begin{proof}
Let $T \subset T^3$ be a convex vertical torus with infinite slope containing $L$, and
let $T_0$ be a nearby translate of $T$. Put $T_0$ in standard form and take a 
convex vertical annulus $A$ between $L$ and a Legendrian divide of $T_0$. The
dividing set $\Gamma_A$ consists of some vertical closed curves and some arcs with
endpoints on $L$. After adapting the characteristic foliation of $A$, we can
take a possibly smaller convex annulus $A' \subset A$ with Legendrian boundary 
$\partial A' = L \cup L_0$ such that $\Gamma_{A'}$ consists only of arcs, then $L$ is a 
stabilisation of $L_0$.
\end{proof}

\begin{prop} \label{trasporto}
If $p_1$, $p_2 \in \Gamma_{T^2}$ are points belonging to dividing curves bounding a 
negative (positive) region $C$, then a positive (negative) stabilisation $L_1$
of $S^1 \times \{ p_1 \}$ is isotopic to a positive (negative) stabilisation $L_2$
of $S^1 \times \{ p_2 \}$.
\end{prop}
\begin{proof}
Let $(T^2 \times [0, \frac 12 ], \xi)$ be a basic slice with 
standard boundary and boundary slopes $s_0=0$ and $s_{\frac 12 }= \infty$ contact 
embedded in
$(T^3, \xi_n)$ so that $S^1 \times \{ p_1 \}$ is a Legendrian divide of 
$T^2 \times \{ \frac 12 \}$ and $T^2 \times \{ 0 \} \subset C \times S^1$. Since the relative
Euler class of $\xi$ is $e(\xi)= \pm \binom{1}{1}$, the sign of the basic slice
$(T^2 \times [0, \frac 12 ], \xi)$ is equal to the sign of the simply connected region
in a convex horizontal annulus of $T^2 \times [0, \frac 12 ]$. It follows that, if 
$C$ is a negative region, then $(T^2 \times [0, \frac 12 ], \xi)$ is a positive basic 
slice and vice versa.

Make the rulings of 
$T^2 \times \{ 0 \}$ vertical, and consider a convex vertical annulus $A$ between
$\{p_1 \} \times S^1 \subset T^2 \times \{ \frac 12 \}$ and a vertical ruling of $T^2 \times \{ 0 \}$. 
The dividing set of $A$ 
consists of a single dividing arc with both endpoints on $T^2 \times \{ 0 \}$, 
and the simply connected region of $A \setminus \Gamma_A$ is positive (negative). Then by 
\cite{etnyre:0}, Lemma 2.20, a vertical Legendrian ruling of $T^2 \times \{ 0 \}$ is a 
positive (negative)
stabilisation of $\{ p_1 \} \times S^1$. From the well definedness of  
stabilisation,
it follows that $L_1$ is contact isotopic to a vertical Legendrian ruling of
$T^2 \times \{ 0 \}$. We can repeat the same argument with a basic slice 
$(T^2 \times [0, \frac 12], \xi)$ with the same sign and the same boundary slopes so 
that  $\{p_2 \} \times S^1$ is a Legendrian divide of $T^2 \times \{ \frac 12 \}$. From this we
 conclude that $L_2$ is contact isotopic to a vertical Legendrian ruling of 
$T^2 \times \{ 0 \}$, then $L_1$ and  $L_2$ are Legendrian isotopic.
\end{proof}
\begin{figure}\centering
\includegraphics[width=10cm]{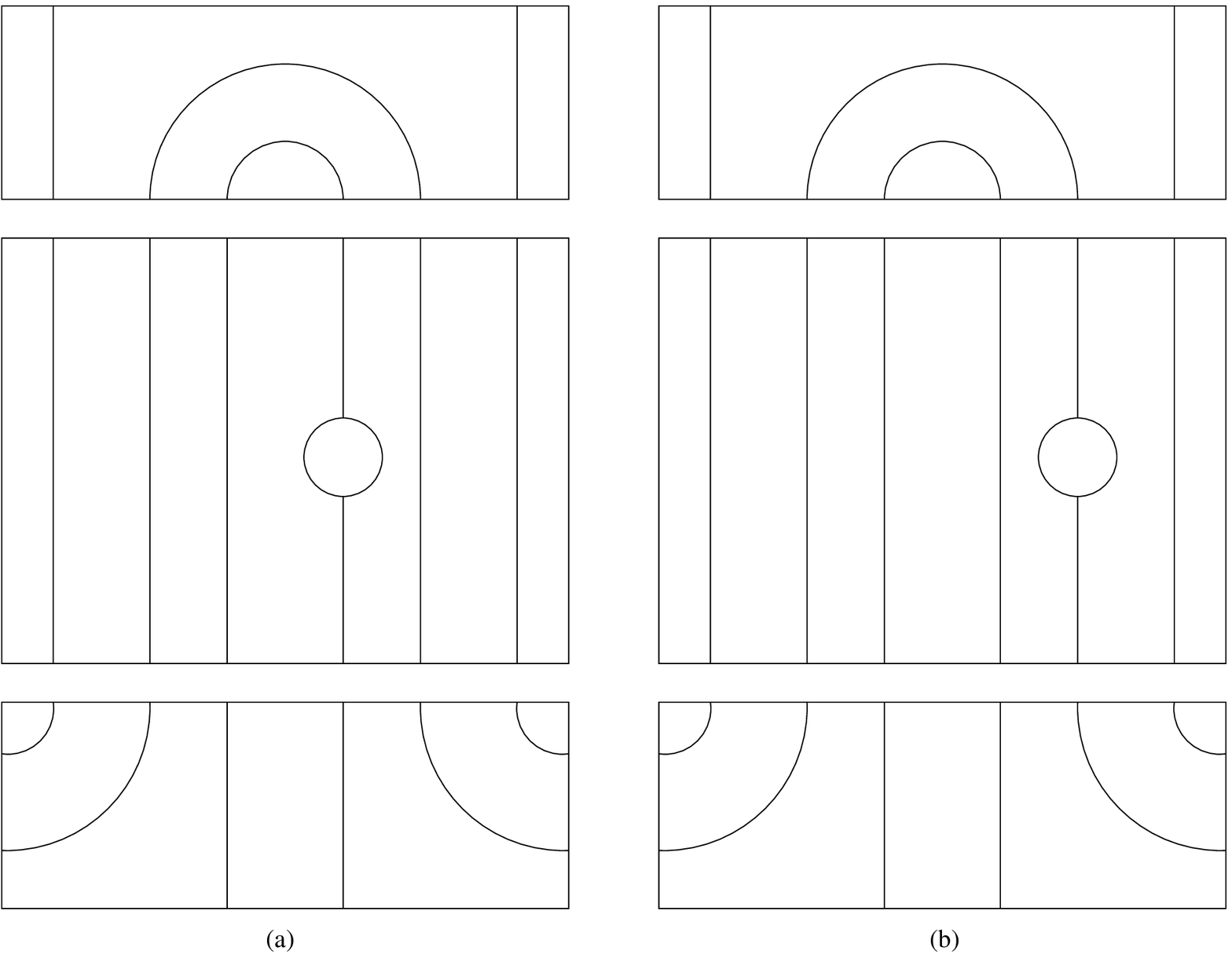}
\caption{The dividing set on a horizontal annulus in $T^2 \times [-1,2] \setminus U_1$ (a) and
in $T^2 \times [-1,2] \setminus U_2$ (b).}
\label{figura}
\end{figure}
   
\begin{prop}
Let $L_1$ and $L_2$ be vertical Legendrian curves in $(T^3, \xi_n)$ with
$\hbox{tb}(L_1)= \hbox{tb}(L_2) <-1$ and $\hbox{r}(L_1)= \hbox{r}(L_2)$. If
$|\hbox{r}(L_i)|<|\hbox{tb}(L_i)|$, then $L_1$ and $L_2$ are Legendrian isotopic. 
\end{prop} 
\begin{proof}
By Proposition \ref{stabilizzazione:t=0} and Proposition 
\ref{massimali-modelli:t=0}, $L_i$ is isotopic to a multiple stabilisation of 
a Legendrian fibre $S^1 \times \{ p_i \}$ with $p_i \in \gamma_i \subset \Gamma_{T^2}$ for $i=1,2$. Because of 
the inequality between the Thurston-Bennequin invariant and the rotation 
number, $L_1$ is obtained by 
applying simple stabilisations of different signs. In particular, being a 
positive stabilisation, $L_1$ is Legendrian isotopic to a multiple 
stabilisation of $S^1 \times \{ p_1' \}$, with $p_1'$ belonging to a dividing curve
$\gamma_1'$ bounding a negative region with $\gamma_1$. In the same way, being a negative 
stabilisation, $L_1$ is Legendrian isotopic to a multiple
stabilisation of $S^1 \times \{ p_1'' \}$, with $p_1''$ belonging to a dividing curve
$\gamma_1''$ bounding a positive region with $\gamma_1$. Going on with the same argument,
we conclude that $L_1$ and $L_2$ are multiple stabilisations of the same 
Legendrian fibres, then  they are Legendrian isotopic because their classical 
invariants agree.
\end{proof}

\begin{prop}\label{fine}
Let $L_1$ and $L_2$ be (multiple) stabilisations of Legendrian fibres 
$S^1 \times \{ p_1 \}$ and $S^1 \times \{ p_2 \}$ respectively with $\hbox{tb}(L_1)= \hbox{tb}(L_2)$,
$\hbox{r}(L_1)= \hbox{r}(L_2)$, and $| \hbox{tb}(L_i)|=| \hbox{r}(L_i)|$. If 
$L_1$ and $L_2$ are isotopic, then either $p_1$ and $p_2$ belong to the 
same dividing curve, or the dividing curves containing $p_1$ and $p_2$ bound 
a region $C$, and the sign of $C$ is opposite to the sign of $\hbox{r}(L_i)$.
\end{prop}

\begin{proof}
First observe that the statement is non trivial only if $\# \Gamma_{T^2} > 2$ because,
if $\# \Gamma_{T^2} = 2$, the two dividing curves bound both a positive and a negative
region. In the proof we will suppose $\mbox{r}(L_i)>0$ and $C$ negative. The 
other case is completely symmetric.
For $i=1,2$, let $U_i$ be a standard neighbourhood of $S^1 \times \{ p_i \}$ containing 
$L_i$,  and  let $V_i \subset U_i$ be a standard neighbourhood of $L_i$. Let $\varphi_t$ be
a contact isotopy of $T^3$ inducing the Legendrian isotopy from $L_1$ to $L_2$. 
We can assume that $\varphi_1(U_1)=U_2$ and $\varphi_1(V_1)=V_2$.   
Lifting the isotopy between $L_1$ and $L_2$ to a finite cover, if needed, we can
also assume that $\varphi_t$ is constant on some vertical torus $T$. We choose $T$ so 
that it is convex with infinite slope and 
$T^3 \setminus T \cong T^2 \times I$ carries an $I$-invariant nonrotative tight contact structure
$\bar{\xi}$, therefore $\varphi_1$ is a 
contactomorphism between $(T^2 \times I \setminus V_1, \bar{\xi}|_{T^2 \times I \setminus V_1})$ and 
$(T^2 \times I \setminus V_2, \bar{\xi}|_{T^2 \times I \setminus V_2})$ fixed on the boundary of $T^2 \times I$. 

From $(T^2 \times I \setminus V_1, \bar{\xi}|_{T^2 \times I \setminus V_1})$ and $(T^2 \times I \setminus V_2, \bar{\xi}|_{T^2 \times I \setminus 
V_2})$ we will construct two isomorphic contact contact
manifolds $(M, \eta_1)$ and $(M, \eta_2)$ and, assuming that $p_1$ and $p_2$ belong to 
different dividing curves which do not bound a negative region, we will prove 
that $(M, \eta_1)$ is tight and $(M, \eta_2)$ is overtwisted, thus obtaining a 
contradiction. 
Glue $S^1$-invariant nonrotative tight contact structures $(T^2 \times [-1,0], 
\bar{\xi}')$ and $(T^2 \times [1,2], \bar{\xi}'')$ with horizontal annuli as in Figure
\ref{figura} to $(T^2 \times I, \bar{\xi})$ such that the boundary parallel
region in the horizontal annulus of $T^2 \times [1,2]$ is negative,
 and call $(T^2 \times [-1,2], \bar{\xi})$ the contact manifolds obtained. 
Let $\alpha_0'$ be a tight contact structure  on a solid torus $S_0'$ with boundary 
slope $-1$, and let $\alpha_{-1}$ be a tight contact structure on a solid torus 
$S_{-1}$ with boundary slope $n+1$ and relative Euler class $n$.
 For $i=1,2$, let  $(M, \eta_i)$ be the contact manifolds obtained from 
$(T^2 \times [-1,2] \setminus V_i, \bar{\xi}|_{T^2 \times [-1,2] \setminus V_i})$ by 
gluing $(S_0', \alpha_0')$ to $\partial V_i$ by 
the map $A_0= \left ( \begin{matrix} n+1 & 1 \\
                         -1  &  0 \end{matrix} \right )$, and gluing
 $(S_{-1}, \alpha_{-1})$ to $T^2 \times \{ -1 \}$, by the map 
$A_{-1} = \left ( \begin{matrix} n+1 & -1 \\
                          1  &  0 \end{matrix} \right )$.
The following lemmas conclude the proof of the theorem.
\end{proof}

\begin{figure}\centering
\includegraphics[width=5cm]{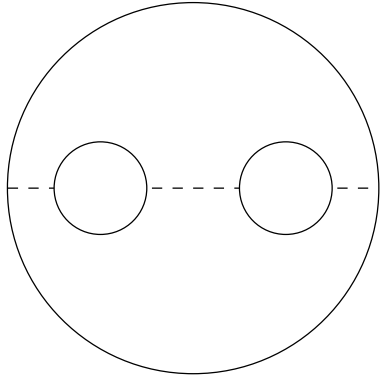}
\caption{The dividing set on a horizontal annulus in $T^2 \times [-1, -1+ \epsilon] \setminus U_1$.}
\label{pantaloni}
\end{figure}

\begin{figure}\centering
\includegraphics[width=12cm]{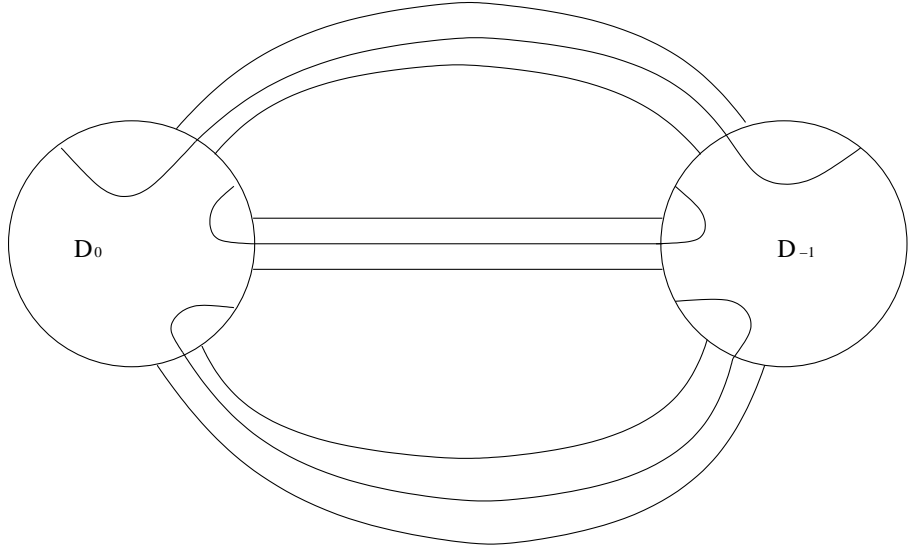}
\caption{The signs of the boundary parallel region of $D_i \setminus \Gamma_{D_i}$ determine
the dividing set of $F'$: an example with $n=2$.}
\label{patching.fig}
\end{figure}

\begin{lemma}
The contact manifold $(M, \eta_1)$ constructed in Proposition \ref{fine} is tight.
\end{lemma}
\begin{proof}
Isotope the product structure on $T^2 \times [-1, 2]$ so that $T^2 \times [-1, -1+ \epsilon]$, for 
$\epsilon >0$ small, becomes an invariant neighbourhood of $T^2 \times \{ -1 \}$.
We can  assume that the standard neighbourhood $U_1$ of $S^1 \times \{ p_1 \}$ is 
contained in $T^2 \times [-1, -1+ \epsilon]$, then the dividing set
on a convex horizontal surface $\Sigma_0 \subset T^2 \times [-1, -1+ \epsilon] \setminus U_1$ looks like in 
Figure \ref{pantaloni}. Call  $M' = M \setminus (T^2 \times (-1+ \epsilon ,2])$ and $\eta' = \eta|_{M'}$.

Let $\alpha_0$ be a tight contact structure on a solid torus $S_0$ with boundary slope
$-(n+1)$ and relative Euler class $n$. We can identify $(S_0', \alpha_0')$ to a 
contact submanifold of $(S_0, \alpha_0)$ so that $A_0$ maps $(S_0 \setminus S_0', \alpha_0|_{S_0 \setminus S_0'})$
onto $(U_1 \setminus V_1, \eta_1|_{U_1 \setminus V_1})$.
The vertical vector field $\frac{\partial}{\partial x}$ on $T^2 \times [-1- \epsilon ,-1] \setminus U_1$ is a 
contact vector field for
$\bar{\xi}|_{T^2 \times [-1- \epsilon ,-1] \setminus U_1}$.
 This vector field extends  to the whole $M'$ to a
vector field $X'$ tangent to the fibres of a Seifert fibration on $M'$.
Since $X'$ is transverse to the meridional discs $D_i$ of $S_i$, it 
is possible to isotope $\alpha_{-1}$ and $\alpha_0$ so that $X'$ becomes  a 
contact vector field for $\eta'$. In fact, $\Gamma_{D_0}$ and $\Gamma_{D_{-1}}$ consist of 
boundary parallel arcs, so by the example discussing bypass attachment in $D^2$
in \cite{gluing} Section 2,
there is no nontrivial bypass attachment when flowing $D_0$ and $D_1$ along $X'$,
therefore by \cite{gluing}, Lemma 2.10, $\alpha_{-1}$ and $\alpha_0$ are isotopic to 
$X'$-invariant contact structures.

The flow of $X'$ defines a covering of $(M', \eta_1')$ whose total space can be 
identified with $F' \times \R$, where $F'$ is a surface transverse to the fibres of 
the Seifert fibration on $M'$. The surface $F'$ can be constructed joining 
$D_{-1}$ to $D_0$ by $n+1$ strips.
A similar construction was used in \cite{etnyre-honda:1} 
(see also \cite{ghiggini-schonenberger}). $F'$ is convex because it is 
transverse to $X'$, and $\Gamma_{F'}$ consists of $n+1$ boundary parallel arcs cutting
out positive region: see Figure \ref{patching.fig}. The vector field $X'$ can 
be extended to a contact vector 
field $X$ on $(M, \eta_1)$ because $(T^2 \times [1+ \epsilon, 2], \eta|_{T^2 \times [1+ \epsilon, 2]})$ is 
$S^1$-invariant. Let $F \times \R$ be the total space of the covering of $M$ defined 
by the flow of $X$, and let $\widehat{\eta}$ be the pull-back of $\eta_1$ to $F \times \R$.
Since $\widehat{\eta}$ is $\R$-invariant, by Giroux's Tightness criterion 
\cite{honda:1}, Theorem 3.5, it is tight
because $\Gamma_F$ contains no homotopically trivial closed curves. In fact, $F \setminus F'$
consists of $n+1$ copies of a horizontal annulus in $T^2 \times [1+ \epsilon, 2]$, whose 
dividing arcs with endpoints on $T^2 \times \{ 1+ \epsilon \}$ cut out negative regions.
Therefore, these arcs join to the dividing arcs on $F'$ to give homotopically
non trivial closed curves in $F$.
\end{proof} 

\begin{lemma}
The contact manifold $(M, \eta_2)$ constructed in Proposition \ref{fine} is 
overtwisted. 
\end{lemma}

\begin{proof}
A convex horizontal annulus in $T^2 \times [-1,2] \setminus U_2$ looks like in Figure 
\ref{figura} (b). In particular, there is a dividing arc with both endpoints
on $T^2 \times \{ -1 \}$ which gives a bypass along $T^2 \times \{ -1 \}$.
 The attachment of this bypass gives a convex torus $T'$ with the slope zero
parallel to $T^2 \times \{ -1 \}$. By the Imbalance Principle, \cite{honda:1}, 
Proposition 3.17,
a vertical convex annulus between a vertical ruling of $T'$ and a Legendrian 
divides of $T^2 \times \{ 2 \}$ produces a bypass along $T'$.  
The attachment of this bypass gives a convex 
torus $T''$ with infinite slope again. The layer between $T^2 \times \{ -1 \}$ and $T''$
has twisting $\pi$, then by \cite{honda:1}, Proposition 4.16, it contains a 
convex torus with slope $\frac{1}{n+1}$ parallel to $T^2 \times \{ -1 \}$. 
A Legendrian divide of this torus bounds an overtwisted disk in $(M, \eta_2)$. 
\end{proof}

\nocite{giroux:1}
\bibliographystyle{amsplain}
\bibliography{contatto}
\end{document}